\documentclass[10pt, a4paper]{amsart}
\input amssym.def
\input amssym.tex
\title[ Ranks
 of the Sylow 2-Subgroups of the 
 Classical Simple Groups     ]
{ $\mbox { Ranks
 of the Sylow 2-Subgroups of the 
 Classical Simple Groups    }$      }

\author{ Mong Lung Lang}

\begin{document}
\baselineskip=11pt


\maketitle
\vspace {-1cm}

\begin{abstract}
Let $S$ be a  2-group. The rank 
of $S$ is the maximal dimension of an  elementary
 abelian subgroup  of $S$ over $\Bbb Z_2$.
 The purpose of this article is to determine
 the  rank of  $S$,
 where $S$ is a Sylow 2-subgroup of the 
 classical simple groups of odd characteristic.
\end{abstract} 

\medskip

\begin{center}
{ 0. Introduction}
\end{center}
A 2-group is called {\em realisable} if it is isomorphic
 to a Sylow 2-subgroup of a finite simple group. It is 
 well known that very few 2-groups are realisable (see [HL], [M]).
 In order 
 to understand such realisable 2-groups, one would like
 to find a set of  invariance of 2-groups which enables
 us to differentiate $\Omega_1$ (the set of realisable
 2-groups) from $\Omega_2$ (the set of non-realisable
 2-groups).
This article studies  the 2-rank 
 (an obvious choice of invariance of 2-groups)
of the 
 classical simple groups of 
 odd characteristic. The results 
 are tabulated in the following table.
 Note that the 2-ranks of 
$PSL_{2n+1}(q)$ and  $PSU_{2n+1}(q)$
 are the same as  
$GL_{2n}(q)$ and  $U_{2n}(q)$
 respectively and the 2-ranks of 
$GL_{2n}(q)$ and  $U_{2n}(q)$
 have been determined in [L].
 Note also that classical   groups of small
 ranks have be determined (Theorem 4.10.5 of
 [GLS3]).
\medskip
{\small \begin{center}
$
\begin{array}{|l|l|} \hline
 \vrule height 8pt width 0pt depth 5pt
\phantom{\Big |}
 \mbox{Group  }  & \mbox{2 Rank }  \\ \hline

\phantom{\Big |}
PSL_{4} (q)\,\,,\,PSU_{4} (q)  & 4\\ \hline
\phantom{\Big |}
PSL_{2n} (q)\, :\, q \equiv 3 \,(4)\,,\,\, n \ge 3 & 2n-2 \\ \hline
\phantom{\Big |}
PSL_{2n}(q)\, :\, q \equiv  1\,(4)\,,\,\, n \ge 3 & 2n-1  \\ \hline
\phantom{\Big |}

PSU_{2n} (q)\, :\, q \equiv 3 \,(4)\,,\,\, n \ge 3 & 2n-1 \\ \hline
\phantom{\Big |}
PSU_{2n}(q)\, :\, q \equiv  1\,(4)\,,\,\, n \ge 3 & 2n-2  \\ \hline
\phantom{\Big |}

PSL_{2n+1}(q)\,,\,\, PSU_{2n+1}(q)            & 2n  \\ \hline
\phantom{\Big |}

PSp_{4}(q)  & 4   \\
\hline
\phantom{\Big |}
PSp_{8}(q) & 6  \\
\hline

\phantom{\Big |}
PSp_{2n}(q)\,:\, n \ne 2, 4 & n +1  \\
\hline
\phantom{\Big |}
P\Omega _{2n+1}(q) & 2n \\
\hline

\phantom{\Big |}
P\Omega _{2n}(\eta, q)\,:\
 q^n \equiv -\eta\,(4)
 & 2n-2 \\
\hline

\phantom{\Big |}
P\Omega _{2n}(\eta, q)\,:\
 q^n \equiv \eta\,(4)\,,\,\,
 n \ge 4
 & 2n-2\\
\hline

\end{array}$
\end{center}}

\section {Projective  Symplectic Groups}
\smallskip
\noindent 
Let $X$ be a 2-group. The  wreath product of $X$ and $\Bbb Z_2$
 is given as follows :
 $$  w(X) =  w_1(X)  =  \left < 
\left (\begin{array}
{cc}
X&0\\
0&1\\
\end{array}\right ),
\left (\begin{array}
{cc}
1&0\\
0&X\\
\end{array}\right )
,
J= \left (\begin{array}
{cc}
0&1\\
1&0\\
\end{array}\right )
\right >.$$
 Let $z_0\in Z(X)-\{1\}$.
For our convenience, we use the following 
notations :
{\small $$ 
\overline X = X/\left <z_0\right >,
 Z =  \left <\mbox{   diag}\,(z_0,z_0)\right >,
\mbox{   diag}\,(X, X)  = B,
\overline B= B/Z,
\overline {w(X)} = w(X)/Z.
\eqno(1.1)$$}
 Note that 
 $$
 B =  \left< 
\left
 (\begin{array}
{cc}
 t &0\\
0&1\\
\end{array}\right )\,:\, t \in X\right >
 \rtimes 
\left< 
\left
 (\begin{array}
{cc}
t&0\\
0&t\\
\end{array}\right )
\,:\, t \in X \right >  = B_1 \rtimes
 B_2.
\eqno (1.2) $$
$$\overline B 
=\overline B_1\rtimes  \overline B_2
\cong B_1\rtimes \overline B_2
\cong X\rtimes \overline X
\,,\,\,
\overline {w(X)}
=\overline B\rtimes \overline {\left<J\right >
} \cong \overline B \rtimes \left <
 J\right >.
\eqno (1.3)$$

\smallskip
\noindent 
{\bf Remark 1.0.} {\bf One should take note that  $\overline X
 =X/\left <z_0\right >$ and the others, such as 
$\overline B_1, \overline B_2,\cdots $  are defined as
$B_1/Z$, $B_2/Z, \cdots $. To be more precise, let $K$ be a 
 group where members of $K$ are $r\times r$ matrices
 over $X$ and let $z_0
 \in Z(X)$ be fixed, the bar notation $\overline K$ is defined
 to be $KZ/Z$,
where $Z =\left <diag\,(z_0,z_0, \cdots , z_0)\right
 >$ $(r$ of them$)$. }

\smallskip
\noindent 
 {\bf Lemma 1.1.}
 {\em Let $X$ be given as 
 in the above. $r_2(\overline {w(X)} )\ge  r_2(\overline X)
 +2.$
Let $\overline E$ be elementary abelian.
 Suppose that $\overline E$ is not a subgroup of 
 $\overline B$. Then $r_2(\overline E)
 \le r_2(\overline X) +2$.
}

\smallskip
\noindent 
 {\em Proof.} Let $ V_0  =  \left< 
\left
 (\begin{array}
{cc}
z_0&0\\
0&1\\
\end{array}\right ),
\left
 (\begin{array}
{cc}
t&0\\
0&t\\
\end{array}\right ),
\left (\begin{array}
{cc}
0&1\\
1&0\\
\end{array}\right )
\,:\, t \in X \right > \subseteq w(X).$ It is clear
 that $r_2(\overline {w(X)})
 \ge r_2(\overline V_0) = 2+r_2(\overline X)$.
This completes the proof of the first
 part of our  lemma. 

\medskip
\noindent
Since $\overline E$ is not a subgroup
 of $\overline B$, 
one has $\overline E = (\overline E\cap
\overline B )\times \overline R$, where $\overline R$
 is not the identity group. Applying Lemma 2.3 of [L],
$\overline R$ is of dimension 1.
 Hence 
 $$r_2(\overline E) = r_2(\overline E \cap \overline B) +1.
\eqno(1.4)$$
Further, 
$ E$ 
 (the preimage of $\overline  E$ in
 $w(X)$) possesses an element
 of the form 
$ \sigma  =\left (\begin{array}
{cc}
0&x_0\\
y_0&0\\
\end{array}\right )\in E -B $
and $\overline E
 = (\overline E\cap \overline B)\times \left <\overline \sigma
 \right >$.
Let
 $\tau=\left (\begin{array}
{cc}
x& 0\\
0&y\\
\end{array}\right )\in E \cap B.$
 Since $\overline E$ is abelian,
 $[\sigma , \tau] \in Z.$
 It follows that 
$x= x_0yx_0^{-1} \mbox{ modulo } \left <z_0\right >.$
 Hence $x \in X$ is uniquely determined by 
$y\in 
X $
 modulo $\left < z_0\right > $.
 As a consequence, the dimension of $ \overline  E \cap \overline B$
 is no more than the rank  of $\overline X +1$
 (1 comes from the fact that $x$ determines $y$ modulo
$  \left <z_0\right >
$, a group of rank  1). Equivalently,
 $r_2( \overline  E \cap \overline B)
 \le   r_2(\overline {X}) +1.
$ The inequality (1.4) now becomes
 $r_2(\overline E) \le 2 +r_2(\overline X)$.\qed

\smallskip
\noindent {\bf Lemma 1.2.} 
 {\em 
Let $Q$ be generalised quaternion
 and let $\overline {w(Q)} = w(Q)/diag(-1,-1).$  Then
$r_2(\overline {w(Q)}) = 4$.
 Let $\overline E$ be
 elementary abelian of dimension $4$. Then $\overline E$
 possesses an element of the form
$ \overline {\left (\begin{array}
{cc}
0&x\\
y&0\\
\end{array}\right ) },$  where $x,y \in Q, xy=yx = \pm 1 .$
}

\smallskip

\noindent {\em Proof.} Let $X=Q$, $z_0 = -1$.
 Applying Lemma 1.1, equation (1.3) and Lemma 2.3 of [L], we have
$$4 = 2+r_2(\overline X) 
 \le  r_2(\overline {w(Q)}   ) 
=r_2(\overline B\rtimes \left<J\right >)
 \le r_2(\overline  B) +1
 \le r_2( X) +r_2(\overline X)+1=4.$$
Since $r_2(\overline B) \le r_2(X) + r_2(\overline
 X) = 1+ 2=3$, $\overline E $ is not a subgroup
 of $\overline B$. Hence $\overline E \cap (\overline
 {w(Q)} - \overline B) \ne \emptyset.$
 Hence $\overline E$ must 
 possess an element of the form
$ \overline {\left (\begin{array}
{cc}
0&x\\
y&0\\
\end{array}\right ) },$  where $x , y\in Q, xy=yx=\pm 1 .$
This completes the proof of the lemma.
\qed

\smallskip
\noindent {\bf Lemma 1.3.} 
 {\em Let $Q$ be generalised quaternion
 and let $\overline {w_2(Q)} = w_2(Q)/diag(-1,-1,
$
$-1,-1).$
Then
$r_2(\overline {w_2(Q)}) = 6$,
$r_2(w(Q)\rtimes \overline {w(Q)})=5.$  In particular,
 if  $\overline E$ is
 elementary abelian of dimension $6$,
 then $\overline E$
 possesses an element of the form
$ \overline {\left (\begin{array}
{cc}
0&x\\
y&0\\
\end{array}\right ) },$  where $x,y \in w(Q), xy=yx = \pm 1 .$

}

\smallskip

\noindent {\em Proof.}
 Let $X= w(Q)$, $z_0 =$ diag$\,(-1,-1)$.
 Then $w_2(Q) = w(X)$. By Lemmas 1.1 and 1.2,
$r_2(\overline {w_2(Q)}) \ge 2+r_2(\overline X) =  6$.
Direct calculation shows that 
 the rank of $ \overline B \cong X\rtimes \overline X
 \cong w(Q)\rtimes \overline {w(Q)}$ is 5.
Let $\overline E$ be an elementary
 abelian subgroup of 
 $\overline {w(X)}$
of dimension $r_2(\overline{w(X)}) \ge 6$.
 It follows from the above that 
 $\overline E$ is not a subgroup of 
 $w(Q)\rtimes \overline {w(Q)}$.
 This implies that $\overline E \cap (
\overline {w_2(Q)} -\overline B) \ne \emptyset$.
 Hence $\overline E$ must 
 possess an element of the form
$ \overline {\left (\begin{array}
{cc}
0&x\\
y&0\\
\end{array}\right ) },$  where $x , y\in w(Q), xy=yx=\pm 1 .$
By Lemmas 1.1 and 1.2, 
 $$r_2(\overline E) \le 
2 +r_2(\overline X) = 2+4=
6 .\eqno(1.5)
$$ 
This completes the proof of the lemma.\qed

\smallskip

\noindent {\bf Lemma 1.4.}
 {\em Let $Q$ be generalised quaternion and let
 $\overline {w_n(Q)} = w_n(Q)/Z$,
 where $z_0 = diag\,(-1, \cdots,-1)\,\,(2^{n-1}\,\,of\,\,them)$.
 Then $r_2(\overline {w_n(Q)}) \ge 2^n+1$.}

\smallskip
\noindent {\em Proof.}
 Let $\left <i,j\right >$ be
 the quaternion subgroup of order 8 of $Q$.
 Let $E_m$ be the matrix obtained from 
 $I_{2^n}$ by replacing the $(m,m)$-entry by
 $-1$ and let 
$$
 M = \left < E_m, \mbox{diag}\,(i,i,\cdots , i),
 \mbox{diag}\,(j,j,\cdots , j), m = 1,2, \cdots, 2^n
\right > \subseteq w_n(Q).$$
 Then $\overline M$ is elementary abelian
 of dimension  $2^n+1$.
 This completes the proof of the lemma. \qed
 
\smallskip

\noindent {\bf Lemma 1.5.} 
 {\em Let $Q$ be generalised quaternion
 and let
 $\overline {w_3(Q)} = w_3(Q)/Z$,
 where $z_0 = diag\,(-1, -1,-1,-1)$.  Then
$r_2(\overline {w_3(Q)}) = 9$.
Let $\overline E$ be elementary abelian of  dimension $9$. 
 Then $\overline E \subseteq diag\,\overline
 {(w_2(Q), w_2(Q))}$.
}
\smallskip

\noindent 
{\em Proof.}
 Let $X = w_{2}(Q)$.
 Let  $\overline E$ be elementary abelian
 of dimension $r_2(\overline {w_3(Q)})$.
 Suppose that $\overline E$ is not a 
 subgroup of $\overline B$ (see (1.2) for notation).
 By Lemmas 1.1 and 1.3,
 $r_2(\overline E) \le r_2(\overline
 X) +2 = 8$.
A contradiction (see Lemma 1.4). Hence $\overline E \subseteq
 \overline B = \mbox{diag}\,\overline
 {(w_2(Q), w_2(Q))
}$. 
 It follows that 
 $$r_2(\overline E)
 \le r_2(\overline B)
=r_2( X\rtimes \overline X)
\le r_2( X) + r_2( \overline X)
= 4+6=10.$$
 Suppose that $r_2(\overline E) = 10$.
 Let $\overline E = \overline V \times 
\overline   W$, where $\overline V = \overline E \cap 
\overline B_1$
 (see (1.2) for notation).
It follows that 
 $\overline V$ is of dimension  $4 = r_2(X)$ and that 
 $\overline {W}$ 
 is of dimension  $6 = r_2(\overline X)$.
 Since $\overline B_1 \cong B_1\cong w_2(Q)$ and elementary
 abelian subgroups of dimension 4 of $w_2(Q)$ are subgroups of 
 $w(Q)\times w(Q)$ (Lemma 3.2 of [L]), we have  
 $$\overline V = \overline
 {\mbox{diag}\,(X_1, X_2, I_4)}
\cong\mbox{diag}\,(X_1, X_2, I_4)
,$$
 where 
$X_i$ are elementary abelian subgroups
 of dimension 2 of $w(Q)$.
An easy observation of the
 structure of $ w(Q)$ 
shows that 
  $X_i$ is either 
$$ \left <
\left (\begin{array}
{cc}
-1&0\\
0&1\\
\end{array}\right ),
\left (\begin{array}
{cc}
1&0\\
0&-1\\
\end{array}\right )\right > \mbox{  or
}
\left <
\left (\begin{array}
{cc}
-1&0\\
0&-1\\
\end{array}\right ),
\left (\begin{array}
{cc}
0& g\\
g^{-1}&0\\
\end{array}\right )\right >.\eqno(1.6)$$
Further, since $r_2(\overline W) =6$, we have 
$r_2(\overline W|\overline B_2) = 6$
 (Lemma 2.3 of [L]). Note that 
$\overline B_2 \cong\overline {w_2(Q)} $.
Applying Lemma 1.3, 
$\overline W|\overline B_2$ possesses an element
 of the form
$$\overline {
\left
 (\begin{array}
{cccc}
0& u&0&0\\
v& 0&0&0\\
0& 0&0&u\\
0& 0&v&0\\
\end{array}\right )}\mbox{ where } u,v\in w(Q),
uv=vu = \pm I_2
.\eqno(1.7)$$
 Since $r_2(\overline W|\overline B_2) = 6 = r_2 (\overline
 B_2) $, $\overline W|\overline B_2$ must contain all the 
involutions of the centre
 of $\overline B_2$. Hence
$$\overline {
\left
 (\begin{array}
{cccc}
-I_2& 0&0&0\\
0& I_2&0&0\\
0& 0&-I_2&0\\
0& 0&0&I_2\\
\end{array}\right )} \in \overline W|\overline B_2.\eqno(1.8)$$
It follows that 
$$\overline {
\left
 (\begin{array}
{ccc}
 x &0&0\\
 0&0&u\\
 0&v&0\\
\end{array}\right )}\,\,, \,\,\,\overline {
\left
 (\begin{array}
{ccc}
 w&0&0\\
 0&-I_2&0\\
 0&0&I_2\\
\end{array}\right )} \in \overline W,
\mbox{ where }  
x, w\in w_2(Q).\eqno(1.9)
$$
 Since the above elements are of order 2 in $\overline
 E$, we have
  $$x^2 = I_4, uv=vu=I_2, w^2 = I_4, \mbox{ or } 
 x^2 = -I_4, uv=vu=-I_2, w^2 = I_4.\eqno (1.10A)$$
 Since the elements in (1.9)
 commute with each other in $\overline E$,  we have
$$xwx^{-1}= -w  \,\,\,(\mbox{not } \pm w).\eqno (1.10B)$$
Note that $[\overline V, \overline W]= \overline 1$.
 It follows that 
$x, w \in
 C_{\overline B_1}(\overline V) \cong 
 C_{ B_1}(\mbox{diag}\,(X_1,X_2, I_4))$.
 An easy observation of the centraliser of 
 $ \mbox{diag}\,(X_1,X_2, I_4)    
     $ in $\overline B_1 \cong B_1 \cong
w_2(Q)$ shows that $C_{\overline B_1}(
\overline V
)    \cong  C_{ B_1}(\mbox{diag}\,(X_1,X_2, I_4))$
 does not possess elements 
 $x,w$ satisfy
$$w^2 = I_4, x^2 = \pm I_4, xwx^{-1}= -w.$$
  A contradiction (see remark).
 Hence $r_2(\overline E) \le 9.$
It follows by Lemma 1.4
 that   $r_2(\overline E) =9$. \qed 

\smallskip
\noindent {\bf Remark.} The centraliser
 of {\small  $\left <
\left (\begin{array}
{cc}
-1&0\\
0&1\\
\end{array}\right ),
\left (\begin{array}
{cc}
1&0\\
0&-1\\
\end{array}\right )\right >$ }
in $w(Q)$ is $Q\times Q$. 
Further, 
$$C_{w(Q)}\left(
\left <
\left (\begin{array}
{cc}
-1&0\\
0&-1\\
\end{array}\right ),
\left (\begin{array}
{cc}
0& g\\
g^{-1}&0\\
\end{array}\right )\right >\right ) =
\left <
\left (\begin{array}
{cc}
 a&0\\
0& g^{-1} a g \\
\end{array}\right ),
\left (\begin{array}
{cc}
0& g\\
g^{-1}&0\\
\end{array}\right )\right >.$$
Direct calculation shows that members of
$w_2(Q) - w(Q)
 \times w(Q)$, which take 
 the form $
\left (\begin{array}
{cc}
0&a_{12}\\
a_{21}&0\\
\end{array}\right )$
 do not centralise diag$\,(X_1,X_2)$.
 It follows that 
$$C_{w_2(X)}(X_1\times X_2)
 \cong  C_{w(X)}(X_1)\times   C_{w(X)}(X_2)   .\eqno(1.11)$$
This enables us to determine the centraliser 
 $ C_{\overline B_1}(\overline V)\cong
C_{B_1}(\mbox{diag}\,(X_1,X_2, I_4))\cong 
C_{w_2(Q)}(\mbox{diag}\,(X_1,X_2))$.
 As a consequence, 
if $x,w \in C_{\overline B_1}(\overline V)$
 gives the relation $w^2 = I_4, x^2 = \pm I_4, xwx^{-1}=\pm w_4$,
 then 
 $xwx^{-1}= w$ all the time !
Note that (1.11) is a special case of the 
 following :
 Let $X$ be a $2$-group and let 
 $E_i$ $(1\le i \le 2^n)$ be subgroups of $X$. Then
$$C_{w_n(X)}(\mbox{diag}\,
 (E_1, \cdots , E_i,\cdots, E_{2^n} ))
 =\mbox{ diag}\,(\cdots, C_X(E_i), \cdots).\eqno(1.12)$$

\smallskip
\noindent 
 {\bf Lemma 1.6.}
 {\em 
 Suppose $n\ge 3$. Then $r_2(\overline {w_n(Q)}) = 2^n+1$.
 Further, 
 elementary abelian subgroups $\overline 
 E$ of dimension
 $\overline {w_n(Q)}$ are subgroups of 
$$ \mbox{diag}
\,(\overline {w_n(Q),  w_n(Q))}.$$}

\smallskip
\noindent {\em Proof.} In the case $n=3$,
  by Lemma 1.5, we have
$r_2(\overline {w_3(Q)}) = 2^3+1$ and 
$overline E \subseteq \mbox{diag}
\,(\overline {w_2(Q),  w_2(Q))}$,
Suppose our lemma holds for $n$.
In the $n+1$ case, 
let  $X = w_n(Q)$ and let $\overline E$ be
 elementary abelian of dimension
 $r_2(\overline {w_{n+1}(Q)})$.
 By Lemma 1.4,  $r_2(\overline E)\ge 2^{n+1} +1$.
 If $\overline E$ is not a subgroup of $\overline
B$ (see (1.2) for notation),
then $r_2(\overline E)\le  r_2(\overline X)
 +2$ (Lemma 1.1).  It follows that 
$$ r_2(\overline X)+2 \ge  r_2(\overline
 E) +2 \ge 2^{n+1}+1.$$
By inductive hypothesis, $ r_2(\overline
 X) =   2^n+1$. Hence
 $2^n+3 \ge 2^{n+1}+1.$
 A contradiction. Hence $\overline E$ is a subgroup of 
 $\overline B = \mbox{diag}
\,(\overline {w_n(Q),  w_n(Q))}$.
Applying (1.2) and Lemma 2.3 of [L],
  $$r_2(\overline E) = r_2(\overline B)
 \le r_2(X) + r_2(\overline X) =
 2^n+2^n+1 = 2^{n+1} +1 .$$
 Hence $r_2(\overline E)
 \le r_2(\overline B) \le 2^{n+1} +1$.
By Lemma 1.4, we have
 $r_2(\overline E) = 2^{n+1} +1$.
  \qed

\medskip
\noindent 
Let $n = 2^{m_1} + 2^{m_2} +\cdots + 2^{m_u}$
 ($m_1 <m_2 <\cdots) $
 be the 2-adic representation of $n$.
 Then
 $$ S =\,\mbox{diag}\,(w_{m_1}(Q),
w_{m_2}(Q),\cdots ,  w_{m_u}(Q))$$
is isomorphic to a Sylow 2-subgroup of 
 $Sp_{2n}(q)$ ($Q$ is generalised
 quaternion). Let
 $Z = \mbox{diag}\,(-1, -1, \cdots , -1)$ ($n$ of them).
 Then $\overline S = S/Z$ is 
isomorphic to a Sylow 2-subgroup of 
 $PSp_{2n}(q)$ (Theorem 6 of Wong [W]).
 Similar to (1.1)-(1.3), 
 we have 
 $$\overline S\cong
  ( w_{m_2}(Q) \times 
w_{m_3}(Q) \times \cdots \times w_{m_u}(Q))
\rtimes \overline {w_{m_1}(Q) }. \eqno(1.13) $$
 To be more precise, one has
 $$
   \overline S = \overline B_1
\rtimes \overline B_2, \eqno(1.14) $$
where
 $B_2 =\{ \mbox{diag}\,(t,t,\cdots, t)\,\,(n/2^{m_1} \mbox{ of 
 them})\,:\,            
 t\in w_{m_1}(Q)\} \cong w_{m_1}(Q)$ and $B_1 = $
$\mbox{diag}\,{(  
w_{m_2}(Q) , \cdots , w_{m_u}(Q), I_{2^{m_1}}   )}.$ Note that 
$$\overline B_1 =\mbox{diag}\,\overline{( 
w_{m_2}(Q) , \cdots , w_{m_u}(Q) , I_{2^{m_1}}   )}\cong
\mbox{diag}\,{(
w_{m_2}(Q) , \cdots , w_{m_u}(Q),  I_{2^{m_1}} )}.$$
Applying the proof of Lemma 1.4,
 we have $$
r_2(\overline S) \ge n+1.\eqno(1.15)$$
 Applying Proposition 3.4 of [L] and
 Lemma 1.6, we have
\begin{enumerate}
\item[(i)] if $m_1=0$, then
$r_2(\overline S) \le
 \sum_{i=2}^{u} r_2(w_{m_i}(Q))+r_2(\overline
 Q) =n+1$, 
\item[(ii)] if $m_1=1$, then
$r_2(\overline S) \le
 \sum_{i=2}^{u} r_2(w_{m_i}(Q)) +r_2(\overline
 {w(Q)}) =n+2$,
\item[(iii)] if $m_1=2$, then
$r_2(\overline S) \le
 \sum_{i=2}^{u} r_2(w_{m_i}(Q)) +r_2(\overline
 {w_2(Q)}) =n+ 2$,
\item[(iv)] if $m_1 \ge 3$, then
$r_2(\overline S) \le
 \sum_{i=1}^{u} r_2(w_{m_i}(Q)) 
 = n+1.$
\end{enumerate}

\smallskip
\noindent (i) By (1.15), we have $r_2(\overline S) = 
 n+1$.

\smallskip
\noindent 
 (ii) The case  $n=2$ is covered by Lemma 1.2. We shall
 assume that $n>2$. 
We have $B_2 \cong w(Q)$ (see (1.14) for notation).
 By Lemma 1.2,
 $r_2(\overline B_2) =4$.
Suppose that  $r_2(\overline S)=n+2 $.
 Let $\overline E$ be elementary abelian of dimension
 $n+2$. Let $\overline E = \overline V \times \overline W$,
 where $\overline V = \overline E\cap \overline B_1$.
The equality
$r_2(\overline S) = r_2(\overline E) = 
 \sum_{i=2}^{u} r_2(w_{m_i}(Q)) +r_2(\overline
 {w(Q)}) =n+2$ implies that $ V\cong \overline V
$ ($V$ is the preimage of 
 $\overline V$) is of dimension $n-2$ and
 $\overline W$ is of dimension 4. 
 Since $r_2(w(Q)) =2$, applying Lemma 3.2 of [L],
 we have 
$$ V 
 = \mbox{diag}\,
(X_1, X_2, \cdots , X_{n/2-1}, I_2),\eqno (1.16)$$
 where $X_i \subseteq w(Q)$, $r_2(X_i) = 2$.
 Applying Lemma 2.3 of [L], $r_2(\overline W)
 = r_2(\overline W|\overline B_2) = 4$.
 As a consequence, $ D= 
 \left < \mbox{diag}\,
(X_1, I_2,  \cdots , I_2, I_2), \overline W|\overline 
 B_2 \,\right >$
 is elementary abelian of dimension 6.
 Note that $D$ is isomorphic to an elementary
 abelian subgroup of $w(Q)\rtimes \overline {w(Q)}$.
 This is a contradiction (see Lemma 1.3).
 Hence $r_2(\overline S) = n+1$.

\smallskip
\noindent 
 (iii)
 The case  $n=4$ is covered by Lemma 1.2. We shall
 assume that $n>4$. 
We have $B_2 \cong w_2(Q)$. By Lemma 1.2,
 $r_2(\overline B_2) =6$.
Suppose that  $r_2(\overline S)=n+2 $.
 Let $\overline E$ be elementary abelian of dimension
 $n+2$. Let $\overline E = \overline V \times \overline W$,
 where $V = E\cap B_1$.
The equality
$r_2(\overline S) = r_2(\overline E) = 
 \sum_{i=2}^{u} r_2(w_{m_i}(Q)) +r_2(\overline
 {w_2(Q)}) =n+2$ implies that $ V\cong \overline V
$ is of dimension $n-4$ and
 $\overline W$ is of dimension 6. 
 Since $r_2(w(Q)) =2$, applying Lemma 3.2 of [L],
 we have 
$$ V 
 = \mbox{diag}\,
(X_1, X_2, \cdots , X_{n/2-1}, I_4),\eqno (1.17)$$
 where $X_i \subseteq w(Q)$, $r_2(X_i) = 2$.
 Applying Lemma 2.3 of [L], $r_2(\overline W)
 = r_2(\overline W|\overline B_2) = 6$.
 As a consequence, $ D= 
 \left < \mbox{diag}\,
(X_1, X_2, I_4,  \cdots , I_4, I_4), \overline W|\overline 
 B_2 \,\right >$
 is elementary abelian of dimension 10.
 Note that $ D$ is isomorphic to an elementary
 abelian subgroup of $w_2(Q)\rtimes \overline {w_2(Q)}$.
 This is a contradiction (see Lemma 1.5).
Hence 
 $r_2(\overline S)=n+1$.

\smallskip
\noindent (iv) By (1.15), we have $r_2(\overline S) = 
 n+1$.

\smallskip
\noindent In summary, we have the following result :

\smallskip
\noindent
{\bf Proposition 1.7.}
{\em Let $S$ be a Sylow $2$-subgroup of 
 $PSp_{2n}(q)$. Then 
\begin{enumerate}
\item[(i)] 
$r_2(S) = 4$ if $n=2$, $r_2(S) = 6$ if $n=4$,
\item[(ii)]  $r_2(S) = n+1$ if $n\ne  2,4$.
\end{enumerate}}

\section {Projective Special Linear and
 Unitary
 Groups  }

\noindent 
\subsection {Projective Special Linear 
 and Unitary Groups I} The main purpose of this
 section is to determine the 2-rank of 
 $PSL_{2n}(q)$, where $q\equiv 3$ (mod 4) and
 the 2-rank of 
 $PSU_{2n}(q)$ where $q\equiv 1$ (mod 4).

\smallskip
\noindent {\bf Lemma 2.1.} 
 {\em Let $K$ be a $2$-group and let $z_0\in Z(K)^{\times},
\overline K = K/\left <z_0\right >$.
Suppose that $r_2(K) = r_2(\overline K) =2$. Then
$r_2(\overline {w(K)}) = 4$.}

\smallskip

\noindent {\em Proof.}  Let $X=K$.
By Lemma 1.1
$r_2(\overline {w(X)}) \ge r_2(\overline X) + 2 \ge 4$.
 Let $\overline E$ be elementary abelian
 of dimension  $r_2(\overline {w(X)})$. Suppose that 
 $\overline E$ is not a subgroup of $\overline B$
 (see (1.1)-(1.3) for notation).
 By Lemma 1.1, $r_2(\overline E) \le 4 $.
Suppose that $\overline E \subseteq \overline B$.
By Lemma 2.3 of [L],
$$r_2(\overline E) =r_2(\overline B)
 =r_2(X\rtimes \overline X) \le 4.$$
 This completes the proof of the lemma. \qed

\smallskip
\noindent {\bf Lemma 2.2.}
 {\em  Let $K$ be a $2$-group. $z_0\in Z(K)^{\times},
\overline {w_n(K)}= w_n(K)/diag\,(z_0,\cdots z_0)$
$
 (2^n\,\,of\,\,them)$.
Suppose that $r_2(K)  =2$. Then
$r_2(\overline {w_n(K)}) \ge 2^{n+1}-1$.}

 \smallskip
\noindent {\em Proof.}
Since $ r_2(K) =2$, $r_2(w_n(K)) = 2^{n+1}$
 (Proposition 3.4 of [L]).
 It follows that $r_2(\overline {w_n(K)})
 \ge 
r_2({w_n(K)})-1 = 
2^{n+1}
-1$.\qed

\smallskip
\noindent {\bf Lemma 2.3.}
 {\em 
 Let $K$ be a $2$-group. $z_0\in Z(K)^{\times},$
$\overline {w_2(K)}=w_2(K)/diag\,(z_0,z_0,z_0,z_0)$.
Suppose that $r_2(K) = r_2(\overline K) =2$,
$r_2(\overline {K\times K}) = 3.$
Suppose further that
$C_{w(X)}(D)$ possesses no   $x,$ $w$ such that 
$w^2 = I_2$, $x^2 = diag\,(t,t)I_2$,
$xwx^{-1} =  diag\,(e,e) w$
$(t\in \left<z_0\right >,e\in \left<z_0\right > -\{1\}
)$,
 where
  $D$ is any elementary abelian subgroup of 
 dimension $4$ of $w(K)$.
 Then
$r_2(\overline {w_2(K)}) = 7$.}

\smallskip

\noindent {\em Proof.} 
Let $\overline E$ be an elementary
 abelian subgroup of 
 $\overline {w_2(K)}$
of dimension $r_2(\overline {w_2(K)})$ and let 
 $X =w(K)$.
  Suppose that  $\overline  E$ 
 not a subgroup of $\overline B$ (see (1.1)-(1.3)
 for notation).
 Then 
$r_2(\overline E) \le r_2(\overline
 {X}) + 2 = 6$ (Lemmas 1.1, 2.1).
This contradicts Lemma 2.2.
 Hence $\overline E$ is a subgroup of $\overline 
 B$. By Proposition 3.4 of [L], $r_2(w(K))=4$.
Hence 
 $$r_2(\overline E) 
=r_2(\overline B) \le
 r_2(B_1) + r_2(\overline B_2)
 = r_2(X) + r_2(\overline X)
 = 8.$$
Suppose that $r_2(\overline E)=8$.
Let $\overline E= (\overline E\cap \overline B_1) \times \overline W.$
 It follows that $r_2(\overline E\cap \overline B_1)
 = r_2(\overline W) = 4.$ By Lemma 2.3 of [L],
 $r_2(\overline W|\overline B_2) = 4 >3 = r_2(\overline 
 {K\times K})$. Hence 
$\overline W|\overline B_2$ possesses an element of the form
$$
\overline {
\left
 (\begin{array}
{cccc}
0& u&0&0\\
v& 0&0&0\\
0& 0&0&u\\
0& 0&v&0\\
\end{array}\right )},$$
where $ u,v\in K,
uv=vu = t$ for some $t\in\left < z_0\right >$.
 Further, since  $r_2(\overline W|\overline B_2) =r_2(\overline
 B_2) = r_2(\overline {w(X)})$
 (see Lemma 2.1), $\overline W|\overline B_2$ must
 contain the following  involution of $Z(\overline B_2)$.
$$\overline {
\left
 (\begin{array}
{cccc}
e& 0&0&0\\
0& 1&0&0\\
0& 0&e&0\\
0& 0&0&1\\
\end{array}\right )},$$
where  $e\in\left < z_0\right >$ is of order 2. As a 
 consequence, $\overline W$ must contain the following elements.
$$
\overline {
\left
 (\begin{array}
{ccc}
 x&0&0\\
 0&0&u\\
 0&v&0\\
\end{array}\right )}\,,\,\,
\overline {
\left
 (\begin{array}
{ccc}
w&0&0\\
 0&e&0\\
 0&0&1\\
\end{array}\right )},$$
where  $x,w\in w(X)$, $x^2 = diag\,(t,t)$, $w^2 = I_2$.
 Further, $xwx^{-1} = diag\,(e,e) w$ ($\overline W$ is 
 abelian).
Since $[\overline E \cap B_1, \overline W] = \overline 1$,
 $x, w \in C_{\overline B_1}(\overline E\cap \overline B_1)$,
 where $\overline B_1 \cong w(X)$ and that 
$r_2(\overline E\cap \overline B_1)=4.$
This contradicts our assumption
 ($\overline E\cap \overline B_1\cong \overline
 E \cap B_1$ is of dimension 4). Hence 
$r_2(\overline E) = r_2(\overline {w_2(T)}) \le 7$.
By Lemma 2.2, we conclude that 
$r_2(\overline {w_2(T)}) = 7$.\qed

\smallskip

\smallskip
\noindent {\bf Lemma 2.4.}
 {\em  Let $K$ be a $2$-group.
Suppose that $r_2(K) = r_2(\overline K) =2$,
$r_2(\overline {K\times K}) = 3.$
Suppose that $r_2(K) = r_2(\overline K) =2$,
$r_2(\overline {K\times K}) = 3.$
Suppose further that
$C_{w(X)}(D)$ possesses no   $x,$ $w$ such that 
$w^2 = I_2$, $x^2 = diag\,(t,t)I_2$,
$xwx^{-1} =  diag\,(e,e) w$
$(t\in \left<z_0\right >, e\in \left<z_0\right >  -\{1\}  )$,
 where
  $D$ be an elementary abelian subgroup of 
 dimension $4$ of $w(K)$.
 If $n \ge 2$, then
$r_2(\overline {w_n(K)}) = 2^{n+1}-1$.}

\smallskip

\noindent {\em Proof.}
Our assertion holds for $n = 2$ (Lemma 2.3).
 Suppose that $r_2(\overline {w_{m} (K)}) = 2^{m+1}-1$.
Let $\overline E$ be 
 elementary abelian of maximal rank
 of $\overline {w_{m+1}(K)}$.
 Let $X = w_m(K)$. Then $w_{m+1}(X) = w(  X   )$.
 Suppose that $\overline E$ is 
 not a subgroup of $\overline {B}$.
 Then $r_2(\overline E)
 \le 2+ r_2(\overline X) = 2^{m+1} +1$ (Lemma 1.1).
 This contradicts Lemma 2.2. Hence 
 $\overline E$ is a subgroup of $\overline
 {B}$. Applying Proposition 3.4 of [L]
 and our inductive hypothesis, we have
 $$r_2(\overline E)\le r_2(\overline {B})
 \le r_2(B_1) + r_2(\overline B_2)
= r_2(X) + r_2(\overline X)
 = 2^{m+1}+ 2^{m+1}-1
 =  2^{m+2}-1.$$
By Lemma 2.2, $r_2(\overline E) =  2^{m+2}-1$. \qed
\smallskip

\noindent Let $T$, $R$, $S(T,R,J)$ and $W(TR, 1,J)$ be
 given as in sections 4.1 and 4.2 of [L].

\medskip
\noindent {\bf Lemma 2.5.} {\em
 Let $T\rtimes R$ be a semidirect product of $2$-groups.
 Suppose that $r_2(TR)$
$ >r_2(T)$ and that 
 $
 \sum_{i=2}^u r_2(w_{m_i}(TR, 1,J)) + r_2(\overline
 {w_{m_1}(TR, 1,J))}
$
$=r_2(\overline {W(TR,1, J)}) 
.$ Then  $r_2(\overline {W(TR,1, J)})
 > r_2(\overline {S(T,R, J)})$.}
\smallskip

\noindent {\em Proof.}
Recall first that 
 $\,\,\overline {
\mbox{diag}\,(w_{m_1}(TR,1,J)\,, w_{m_2}(TR,1,J)\,,
\cdots \,, w_{m_u}(TR,1,J))}$
 $
 =\overline {W(TR,1,J)}$.
 It is easy to see that the above can be 
 written as the following semidirect product.
$$\mbox{diag}\,(w_{m_2}(TR,1,J), w_{m_3}(TR,1,J),
\cdots , w_{m_u}(TR,1,J), I_{2^{m_1}})
\rtimes \overline D,$$
 where $D = \left < \mbox{diag}\,
 (t,t,\cdots)\,:\,
 t \in  w_{m_1}(TR, 1,J)\right >.$
Suppose that 
 $r_2(\overline {W(TR,1, J)})
 = r_2(\overline {S(T,R, J)}) =d$.
Let $\overline E \subseteq\overline {S(T,R, J)}
 $ be elementary abelian of dimension $d$.
 It follows  that 
 $\overline E$ is of maximal dimension in
$\overline {W(TR,1, J)}$.
By the assumption of the lemma, we have
$
r_2( \overline
 {\mbox{diag}\,(I_{2^{m_2}}, \cdots , w_{m_i}(TR,1,J),
\cdots , I_{2^{m_u}}, I_{2^{m_1}})\cap E})
 = r_2(w_{m_i}(TR,1,J)).$
Since $r_2(TR)>r_2(T)$, we have $r_2(TR) \ge 2$.
By Lemma 3.2 of [L], the intersection 
$\overline
 {\mbox{diag}\,( w_{m_2}(TR,1,J), I_{2^{m_3}},
\cdots , I_{2^{m_u}}, I_{2^{m_1}})\cap E})$ must take the following 
 form :
$$\overline
 {\mbox{diag}\,(\prod_{i=1}^{2^{m_2}}
 V_i,I_{2^{m_3}},
\cdots , I_{2^{m_u}}, I_{2^{m_1}})},\mbox{ where }
 V_i\subseteq TR,\,
 r_2(V_i) =r_2(TR)\mbox{ for all } i.$$
It follows that 
 the intersection, as well as $\overline E$,  must 
 contain an element of the 
 form $\overline
 { \mbox{diag}\,(tr, 1,1,\cdots, 1)}$ $(t\in T, r\in R)$.
Since $r_2(TR)>r_2(T)$, $r_2(TR) = 
 r_2(V_i)$, we may assume that $r \ne 1$.
 This implies that $E$, as well as $S(T,R,J)$,
 contains an element
 diag$\,(x_1, x_2, \cdots, )$ ($x_i = t_ir_j, t_i
 \in T, r_i\in R)$, where the number  of $i$ such that
 $r_i\ne 1$ is odd.
 This is a contradiction (see remark of section 4.2 of [L]).
Hence 
 $r_2(\overline {W(TR,1, J)})
 > r_2(\overline {S(T,R, J)})$.\qed

\smallskip

\noindent 
 Let $n = 2^{m_1} + 2^{m_2} +\cdots + 2^{m_u}$
 ($m_1 <m_2 <\cdots) $
 be the 2-adic representation of $n$ and 
 let $T$ be generalised quaternion 
 of order $2^{t+1}$, where $ 2^{t+1}||(q^2-1)$,
and let $R$ be given as in section 5.1 of [L].
Then the group $TR$ satisfies the assumption of 
 Lemma 2.4. Further, 
 $$
S(T,R,J) = 
\mbox{diag}\,( w_{m_1}(T,R,J),
 w_{m_2}(T,R,J),  \cdots ,  w_{m_u}(T,R,J))
\rtimes U(R)$$
 is a Sylow 2-subgroup of $SL_{2n}(q)$,
 $ q\equiv 3$ (mod 4) (Theorem 4 of Wong [W]).
Let  $z_0\in T$ be the element of order 
 $2$.
 Then 
 $Z= \left <\mbox{diag}\,(z_0, z_0, \cdots )
 \right > \subseteq 
  S(T,R,J)$ is the centre of  $S(T,R,J)$.
 Further, 
$$S(T,R,J)/Z\eqno(2.1)$$
 is a Sylow 2-subgroup of $PSL_{2n}(q)$,
 $ q\equiv 3$ (mod 4).
One sees easily that 
$$S(T,R, J)/Z  \subseteq W(TR,1,J)/Z
= \overline {  W(TR,1,J)}. \eqno (2.2)  $$
 Since $W(TR,1,J)\cong \prod w_{m_i}(TR)$ is a direct
 product of wreath products ((4.8) of [L]), one may 
 Proposition 3.4 of [L] to  conclude that 
$$r_2(\overline {W(TR,1,J)}) \ge 
r_2({W(TR,1,J)})-1 = 
2n-1.\eqno (2.3A) $$
Note that 
$\overline {  W(TR,1,J)}$ can be written as a semidirect
 product (see the proof of Lemma 2.5).
$$\overline {W(TR,1,J)}
 \cong \left ( \prod_{i=2}^u
 w_{m_i}(TR) \right )
\rtimes { w_{m_1}(TR)}.$$
 Applying Lemma 2.4 and Proposition 3.4 of [L],
Hence $$r_2(\overline {W(TR,1,J)})
 \le \sum r_2(w_{m_i}(TR)) + r_2(\overline{ w_{m_1}(TR)})
\le 2n-1.\eqno (2.3B)$$
 Hence (see $(2.3A)$  and $(2.3B)$)
 $$r_2(\overline {  W(TR,1,J)}) 
 = 2n-1.\eqno (2.4) $$

\smallskip
\noindent {\bf Proposition 2.6.}
 {\em 
Suppose that $q\equiv 3\,\,\,(mod\,\,4)$. Then 
$r_2(PSL_{4}(q)) =  4.$ If $n \ge 3$, then
 $r_2(PSL_{2n}(q)) = 2n-2.$}

\smallskip
\noindent {\em Proof.}
It is well known that $r_2(PSL_4(q))=4.$
 Suppose that $n \ge 3$.
 $\overline {S(T,R,J)}$ is a Sylow 2-subgroup
 of $PSL_{2n}(q)$.
 Since $r_2(SL_{2n}(q)) = 2n-1$ (see [L]), 
 $r_2(PSL_{2n}(q))
$
$ \ge 2n-2$.
 By (2.4), we have
$$2n-2\le r_2(\overline {S(T,R,J)})
\le r_2(\overline {  W(TR,1,J)}) =2n-1.\eqno(2.5)$$
Since $T$ and $R$ (for $PSL_{2n}(q))$
 satisfy the assumption of Lemma 2.5,
we conclude that 
$ r_2(\overline {S(T,R,J)})
< r_2(\overline {  W(TR,1,J)})$.
Hence $2n-2 = r_2(\overline {S(T,R,J)})$.\qed

\smallskip
\noindent Suppose that 
 $q\equiv 1$ (mod 4). Then (2.1)  is a 
 Sylow 2-subgroup of 
 $PSU_{2n}(q)$. It follows that

\smallskip
\noindent {\bf Proposition 2.7.}
 {\em 
Suppose that $q\equiv 1\,\,\,(mod\,\,4)$. Then 
$r_2(PSU_{4}(q)) =  4.$ If $n \ge 3$, then
 $r_2(PSU_{2n}(q)) = 2n-2.$}

\noindent 
\subsection {Projective Special Linear 
 and Unitary Groups II} The main purpose of this
 section is to determine the 2-rank of 
 $PSL_{2n}(q)$ where $q\equiv 1$ (mod 4) and
 the 2-rank of 
 $PSU_{2n}(q)$ where $q\equiv 3$ (mod 4).

\smallskip
\noindent {\bf Lemma 2.8.}
 {\em Let $K$ be a $2$-group. $z_0 \in Z(K)^{\times}$.
Suppose that $r_2(K) = r_2(\overline K) =2$,
$r_2(\overline {K\times K}) = 4.$
 Then 
$r_2(\overline {w_n(K)}) = 2^{n+1}$.}

\smallskip

\noindent {\em Proof.} By our assumption
 and lemma 2.3 of [L], $\overline {K\times K}$
 has an elementary abelian subgroup of dimension 4 generated by
 $\overline g_i$ ($1\le i\le 4$), where
 $g_1 = $ diag$\,(k_1, 1)$,  $g_2 = $ diag$\,(k_2, 1)$,
 $g_3 = $  diag$\,(a, k_3)$,
 $g_4= $ diag$\,(b, k_4)$, where $k_i, a, b \in K$.
Let $V$ be the subgroup of $w_n(K)$ generated by
 diag$\,(a,a, \cdots , a, k_3)$,
 diag$\,(b,b, \cdots , b, k_4)$, and 
 diag$\,(x_1,x_2, \cdots , x_{2^n-1},  1)$,
 where $x_i$ is either $k_1$ or $k_2$. It is clear
 that $\overline V \subseteq \overline 
 {w_n(K)}$ is elementary abelian of dimension
 $2^{n+1}$. Hence 
 $$r_2(\overline {w_n(K)}) \ge r_2(\overline
 V) = 2^{n+1}.\eqno (2.6)$$
We shall now apply induction on $n$.
It is clear that our assertion holds for 
 $n=0$ as $\overline {w_0(K)} =\overline
  K$ is of rank 2. 
 Suppose that $r_2(\overline {w_{m} (K)}) = 2^{m+1}$.
Let $\overline E$ be 
 elementary abelian of dimension
 $r_2(\overline {w_{m+1}(K)})$.
 Let $X = w_m(K)$. Then $w_{m+1}(K) = w(X)$.
 Suppose that $\overline E$ is 
 not a subgroup of $\overline {B}$
 (see (1.1)-(1.3) for notation).
 Then $r_2(\overline E)
 \le 2+ r_2(\overline X) = 2^{m+1} +1$ (Lemma 1.1).
 This contradicts (2.6). Hence 
 $\overline E$ is a subgroup of $\overline
 {B} $. Applying Proposition 3.4 of [L]
 and our inductive hypothesis, we have
 $$r_2(\overline E)=r_2(\overline {B})
 \le r_2(X) + r_2(\overline X)
 = 2^{m+1}+ 2^{m+1}
 =  2^{m+2}.$$
By (2.6),  
$r_2(\overline E) =  2^{m+2}.$ \qed

\medskip

\noindent 
 Let $n = 2^{m_1} + 2^{m_2} +\cdots + 2^{m_u}$
 ($m_1 <m_2 <\cdots) $
 be the 2-adic representation of $n$.
 Let $T = \left <v, w\right >$ be generalised quaternion 
 of order $2^{t+1}$ $(o(v) = 2^t)$, where $ 2^{t+1}||(q^2-1)$,
 and let $R = \left <e\right >
$, $E_0$, $R_0$ be given as in section 5.2 of [L].
 Then $TR$ satisfies the assumption of Lemma 2.8.
 Further,
  $$
S(T,R,J) =\mbox{diag}\,
( w_{m_1}(T,R,J),
 w_{m_2}(T,R,J), \cdots , w_{m_u}(T,R,J))
\rtimes U(R)$$
 is a Sylow 2-subgroup of $SL_{2n}(q)$,
 $ q\equiv 1$ (mod 4).
Let $2^m = gcd\,(2^{m_1}, q- 1)$ and let
$z_0= (e^2v)^{2^{t-m}}$.
Then the centre
 of $S(T,R,J)$ is generated by
$$ z = \mbox{diag}\,
 (z_0, z_0, \cdots , z_0), \,\,o(z_0)=2^m. \eqno (2.7)
$$ 
By Theorem 6 of Wong [W], 
$$S(T,R,J)/\left <z\right > \eqno(2.8)$$
 is a Sylow 2-subgroup of $PSL_{2n}(q)$,
 $ q\equiv 3$ (mod 4).
One sees easily that 
$$S(T,R, J)/\left <z\right > 
 \subseteq W(TR,1,J)/\left <z\right >
= \overline {  W(TR,1,J)}. \eqno (2.9)  $$
 Since $W(TR,1,J)\cong w_{m_i}(TR)$ is a direct
 product of wreath products ((4.8) of [L]) and
 $K=TR$ satisfies the assumption of Lemma 2.8, similar
 to Lemma 2.8, one may  conclude that 
$$r_2(\overline {W(TR,1,J)})= 2n.\eqno (2.10) $$

\smallskip
\noindent {\bf Proposition 2.9.}
 {\em 
Suppose that $q\equiv 1\,\,\,(mod\,\,4)$. Then 
 $r_2(PSL_{4}(q)) =4 .$ If $n \ge 3$, then
 $r_2(PSL_{2n}(q)) = 2n-1.$}

\smallskip
\noindent {\em Proof.}
 Note that it is well known that 
 $r_2(PSL_{4}(q)) =4 .$ Suppose that $n\ge 3$.
Let
 $$ V
 = \left
 < S(E_0, R_0, 1) \times U(R_0),
\mbox{diag}\,
 (w,w, \cdots , w)
 \right >.$$
 Applying Lemma 4.3 of 
 [L],
we have
 $$
r_2(PSL_{2n}(q)) \ge 
r_2(\overline V) = 2n-1.\eqno(2.11)$$
 $\overline {S(T,R,J)}$ is a Sylow 2-subgroup
 of $PSL_{2n}(q)$. Applying the above
 ((2.10, (2.11)), we have
$$2n-1\le r_2(\overline {S(T,R,J)})
\le r_2(\overline {  W(TR,1,J)}) =2n.\eqno(2.12)$$
Since $T$ and $R$ (for $PSL_{2n}(q))$
 satisfy the assumption of Lemma 2.5,
we conclude that 
$ r_2(\overline {S(T,R,J)})
< r_2(\overline {  W(TR,1,J)})$.
Hence $2n-1 = r_2(\overline {S(T,R,J)})$.\qed

\smallskip
\noindent Suppose that 
 $q\equiv 3$ (mod 4). Then (2.8)  is a 
 Sylow 2-subgroup of 
 $PSU_{2n}(q)$. It follows that

\smallskip
\noindent {\bf Proposition 2.10.}
 {\em 
Suppose that $q\equiv 3\,\,\,(mod\,\,4)$. Then 
$r_2(PSU_{4}(q)) =  4.$ If $n \ge 3$, then
 $r_2(PSU_{2n}(q)) = 2n-1.$}

\section{Orthogonal Commutator Groups $\Omega
 _{2n+1}( q)
 =
 P\Omega_{2n+1}(q)
$} Let $2^{t+1}$ be the greatest power
 of 2 that divides $q^2-1$ and let $T =\left < v,w\right>
$ be a dihedral group of order $2^{t}$,
 where $o(v)= 2^{t-1}, o(w)= 2, wvw= v^{-1}$.
 Further, $R = \left < e\right >$ is a group of 
 order 2 acts on $T$ by
 $eve= v^{-1}, ewe= vw.$
Let 
  $n = 2^{m_1} + 2^{m_2} + \cdots  +2^{m_u}$ 
 be the $2$-adic representation of $n$.
 Then $S(T,R,J)$  is 
 a Sylow 2-subgroup of $\Omega_{2n+1}(q) =
 P\Omega_{2n+1}(q)$
 (see (ii) of Theorem 7 of Wong [W]).
By the results in section 5.4 of [L],
the rank of $S(T,R,J)$ is $2n$.

\section{Orthogonal Commutator Groups $\Omega
 _{2n}(\eta, q)=P\Omega
 _{2n}(\eta, q)   $ where 
 $\eta = \pm 1$, $q^n \equiv -\eta$
 (mod 4)}
Applying Theorem 7 of Wong
 [W],  a Sylow 2-subgroup of 
 $\Omega_{2n}(\eta, q)
 = P\Omega_{2n}(\eta, q)$ is isomorphic 
 to a Sylow 2-subgroup of 
 $O_{2(n-1)} (\eta ',  q)$, where
 $q^{n-1}\equiv \eta '$ (mod 4). Let
 $S$ be a Sylow 2-subgroup of  $O_{2(n-1)} (\eta ',
 q)$, where
 $q^{n-1}\equiv \eta '$ (mod 4).
Applying Theorem 3 of Carter and Fong [CF],
 $S$ is isomorphic to a Sylow
 2-subgroup of $O^+_{2n-1}(q)$.
We shall now describe $S$  as follows : 

\smallskip
\noindent 
Let $D$ be a dihedral group of order $2^{s+1}$,
 where $2^{s+1}$ is the greatest power of 2
 that divides $q^2-1$. Then 
 $D$ is isomorphic to a Sylow 2-subgroup of 
 $O^+_3(q)$. Let $T_{r-1}$ be the wreath product
 of $r-1$ copies of $\Bbb Z_2$ and let
 $S_r$ be the wreath product of $D$ and $T_{r-1}$.
 Then $S_r$ is a Sylow 2-subgroup of 
 $O_{2^r+1}^+(q)$.
 Let $2(n-1) = 2^{m_1} +  2^{m_2}+ \cdots + 2^{m_u}$
be the 2-adic representation of $2(n-1)$. Applying 
 Theorem 2 of Carter and Fong [CF],
 $$S \cong
 S_{m_1} \times  S_{m_2}\times \cdots \times  S_{m_u}.$$
By the results of section 5.5 of [L],
 the rank of $S$ is $2n-2$.

\section{Orthogonal Commutator Groups $P\Omega
 _{2n}(\eta, q)$, where $n$ is even,
 $\eta = \pm 1$, $q^n \equiv \eta$
 (mod 4)} 
\noindent Note first that since $n$ is even and 
 $q^n \equiv \eta$, we have  $\eta =1$.
Let $2^{t+1}$ be the greatest power
 of 2 that divides $q^2-1$ and let
$T$ be the central product of two dihedral
 groups of order $2^{t+1}$ :
$$T = \left
< d= \left (\begin{array}
{cc}
u&0\\
0&u^{-1} \\
\end{array}\right ),
g=\left (\begin{array}
{cc}
u&0\\
0&u \\
\end{array}\right ),
h= \left (\begin{array}
{cc}
0&1\\
1&0 \\
\end{array}\right ),
k=\left (\begin{array}
{cc}
0&w\\
w&0 \\
\end{array}\right )\right>,$$
 where
$o(u) = 2^t$, $o(w)=2$, $wuw=u^{-1}$.
Let
$ R = \left < e, f\right > \cong \Bbb Z_2
 \times \Bbb Z_2$, where 
$$d^e=g^{-1}, g^e=d^{-1}, h^e=gk, k^e=dh,$$
$$d^f=g, g^f= d, h^f=k,k^f=h.$$
Let 
 $n/2 = 2^{m_1} +  2^{m_2}+\cdots +  2^{m_u}$ $(
 n$ even, $m_1 <m_2 <\cdots$) be 
 the 2-adic representation of $n$.
By Theorem 11 of [W], 
 $S(T,R,J)$
  is a Sylow 2-subgroup of $\Omega_{2n}(\eta,
 q)$, where $n$ is even,
 $\eta = \pm 1$, $q^n \equiv \eta$
 (mod 4). 
Let $$z = \mbox{diag}\,
(g^{2^{t-1}},g^{2^{t-1}},
 \cdots ).\eqno (5.1)$$
 Then $\overline {S(T,R, J)}
 = S(T,R, J) /\left < z\right >$ is a Sylow
 2-subgroup of 
 $P\Omega_{2n}(\eta,
 q)$, where $n$ is even,
 $\eta = \pm 1$, $q^n \equiv \eta$
 (mod 4). 
It is easy to see that  
 $$T
 = (\left < d\right>
 \times 
 \left < dg\right>)
\rtimes \left < h,k\right >.$$
 As a consequence, one can show that $r_2(T) =3$,
 $r_2(TR) = 4 = r_2(T) + 1$.
 It is also easy to see that 
 $Z(T) =\left <  g^{2^{t-1}}\right >$ is of order 2 and
 that $r_2(T/Z(T)) = 4$.

\smallskip
\noindent 
 {\bf Lemma 5.1.}
 {\em Let $z_0 =  g^{2^{t-1}}$,
 $\overline {TR} = TR /\left <  z_0\right >$,
$\overline {TR\times TR}
 = (TR\times TR)/\left < diag\,(z_0, z_0)\right >$.
 Then 
$r_2(\overline {TR}) = 4$, 
$r_2(\overline {TR\times TR}) = 7$.
In particular, $r_2(\overline {w(TR)}) = 7.$
}

\smallskip
\noindent 
 {\em Proof.} Let $X = TR$. Then
 $B \cong  TR\times TR$ (see (1.1)-(1.3) for notation).
It is easy to show that $
r_2(\overline X) = r_2(\overline {TR}) = 4$.
 By Lemma 2.3 of [L],
  $$7 =  r_2(B) -1 \le r_2(\overline {B}) \le 8.$$
Suppose that $r_2(\overline {B}) =8$.
 Let $
\overline E \subseteq \overline {B}$
 be elementary abelian of dimension 8
 and let
 $\overline E= \overline V \times \overline W$,
where $\overline V = \overline E \cap \overline B_1$.
 Then both
 $\overline V$ and 
 $\overline W$ are  of dimension 4.
Let
$$\overline V =
\left< 
\overline {\left
 (\begin{array}
{cc}
x_i&0\\
0&1\\
\end{array}\right )}\,:\, 1\le i\le 4 \right >,
\overline W = 
\left< 
\overline {\left
 (\begin{array}
{cc}
u_i&0\\
0&w_i\\
\end{array}\right )} \,:\, 1\le i\le 4\right >.$$
Since $[\overline V, \overline W] =
\overline I_2$, 
$$[x_i, u_j]= 1. \eqno(5.2)$$
Applying Lemma 2.3, $r_2(\overline W|\overline B_2) = 4.$
 Since $r_2(\overline B_2) =4$, 
$\overline W|\overline B_2$ must contain all the involutions
 of $Z(\overline B_2)$. It follows that 
$$\overline {\left
 (\begin{array}
{cc}
g^{2^{t-2}}&0\\
0&  g^{2^{t-2}} \\
\end{array}\right )} \in \overline W|\overline B_2.
\eqno(5.3)$$
 Hence 
$$\overline {\left
 (\begin{array}
{cc}
\sigma &0\\
0&  g^{2^{t-2}} \\
\end{array}\right )} \in \overline W,
 \eqno(5.4)$$
for some $\sigma$.
 Since the above element if of order 2 in $\overline E$, 
 $\sigma^2 = g^{2^{t-1}}$.
 By our results in 
Appendix A, 
$$\sigma \in \{\,
g^{2^{t-2}},\,\,\, d^{2^{t-2}} ,\,\,\,
g^{2^{t-2}}d^mh,\,\,\,
d^{2^{t-2}}g^nk\,\}.\eqno (5.5)$$
 Let 
{\small $A= \left \{ x\,:\,
\overline {\left (\begin{array}
{cc}
x&0\\
0&1\\
\end{array}\right )} \in \overline V \right \}
 \subseteq TR.$}
 Since $A\subseteq TR
  $ is elementary abelian
 of dimension 4 and $T$ is of rank 3, $A$ must
 possess an element (of order 2) of the form
 $tr$, where $t \in T$, $r\in R-\{1\}$.
 This particular element $tr$
must commute with $\sigma$ (see (5.2)).
Since  such elements $tr$ are completely known
 (see b(i) and its remark
 of section 5.6 of [L]) and they do not
 commute with 
$$g^{2^{t-2}}d^mh,\,\,\,
d^{2^{t-2}}g^nk,$$
 (5.5) can be refined into 
$$\sigma \in \{\,
g^{2^{t-2}},\,\,\, d^{2^{t-2}}
\,\}.\eqno (5.6)$$
 Since $[A, \sigma ]=1$ (see (5.2)),
 $A \subseteq C = C_{TR}(\sigma )$. 
This is not possible as such $C$ is 
 of rank at most 3 (see Appendix A).
As a consequence, 
 $r_2(\overline {B}) \ne 8$.
 It follows that $r_2(\overline 
 {(TR\times TR)}) = r_2(\overline {B}) =7$.

\smallskip
\noindent 
Let $\overline E$ be elementary abelian of 
 dimension $r_2(\overline {w(TR)})\ge 7$.
 Suppose that $\overline E $ is not a subgroup of 
$\overline {B}$. By Lemma 1.1,
 $r_2(\overline E)\le 2+r_2(\overline {X})
 = 6 <7$. A contradiction. Hence 
 $\overline E $ is a subgroup of 
$\overline {B}$. This completes the proof
 of the lemma. \qed

\smallskip
\noindent 
Similar to 
Lemmas 1.5 and 1.6, 
which we prove that $r_2(\overline {w_n(Q)}) = 2^n+1$
 under the assumption that 
 $r_2(w_2(Q)) = 4$, $r_2(\overline {w_2(Q)}) = 6$,
$r_2(\overline {w_3(Q)}) = 9$,
one may apply our result
 of Lemma 5.1 ($r_2(\overline {w(TR)}) = 7$), 
$r_2(TR)= r_2(\overline {TR})=4$ to show that 
$r_2(\overline {w_n(TR)}) 
$
$= 2^{n+2}-1$.
Consequently, one has,

\smallskip
\noindent 
 {\bf Lemma 5.2.} Let $T$ and $R$ be given as in the above.
 Then 
 {\em $r_2(\overline {W(TR,1,J)}) = 2n-1$.}

\smallskip
\noindent {\em Proof.}
 Since $W(TR, 1,J) \cong \prod w_{m_i}(TR)$ is a direct product of 
 wreath products, one may apply Proposition 3.4 of [L] 
 to conclude that $r_2(W(TR, 1,J)) = 2n$. 
 Since $r_2(W(TR,1,J))$
 $ = 2n$, we have 
$$r_2(\overline {W(TR,1,J)}) \ge 
r_2( {W(TR,1,J)}) -1=
 2n-1.\eqno(5.7)$$
Recall that
 $$\overline {W(TR,1,J)}
 =\overline B_1 \rtimes \overline B_2
\cong B_1 \rtimes \overline B_2 ,\eqno(5.8)$$
where
$$B_1 =\mbox{ diag}\,(
 w_{m_2}(TR, 1,J),\cdots, 
 w_{m_u}(TR, 1,J), I_{2^{m_1}} )\,,$$$$
 B_2 = 
\left < \mbox{diag}\,
 (t,t,\cdots)\,:\,
 t \in   w_{m_1}(TR, 1,J)    \right > \cong 
 w_{m_1}(TR, 1,J).\eqno(5.9)$$

By  Proposition 3.4 of [L]
 and Lemma 5.1,
 we have
\begin{enumerate}
\item[(i)] if $m_1=0$, then
$r_2(\overline {W(TR,1,J)}) \le
 \sum_{i=2}^{u} r_2(w_{m_i}(TR))+r_2(\overline
 {TR}) = 2n$, 
\item[(ii)] if $m_1\ge 1$, then
$r_2(\overline {W(TR,1,J)}) \le
 \sum_{i=2}^{u} r_2(w_{m_i}(TR)) +r_2(\overline
 {w(TR)}) = 2n-1$.

\end{enumerate}

\smallskip
\noindent 
 (i)  Note first that $B_2 \cong TR$.
Suppose that  $r_2(\overline {W(TR,1,J)})=2n $.
 Let $\overline E$ be elementary abelian of dimension
 $2n$. Let $\overline E = \overline V \times \overline W$,
 where $\overline V = \overline E\cap \overline B_1$.
The equality
$r_2(\overline {W(TR,1,J)}) = r_2(\overline E) = 
 \sum_{i=2}^{u} r_2(w_{m_i}(TR)) +r_2(\overline
 {TR}) =n+2$ implies that $ \overline V
$ is of dimension $2n- 4$ and
 $\overline W$ is of dimension 4. 
 Since $r_2(TR) =4$, applying Lemma 3.2 of [L],
 we have 
$$ \overline V 
 = \mbox{diag}\,
(E_1, E_2, \cdots , E_{n/2-1}, 1),\eqno (5.10)$$
 where $E_i \subseteq TR$ is elementary abelian of 
 dimension 4.
 Applying Lemma 2.3 of [L], $r_2(\overline W)
 = r_2(\overline W|\overline B_2) = 4$.
 As a consequence, $ D= 
 \left < \mbox{diag}\,
(E_1, 1,  \cdots , 1, 1), \overline W|\overline 
 B_2 \,\right >$
 is elementary abelian of dimension 8.
 Note that $D$ is isomorphic to an elementary
 abelian subgroup of $\overline {w(TR)}$.
 This is a contradiction (see Lemma 5.1).
 Hence $r_2(\overline {W(TR,1,J)}) < 2n$. By (5.7),
 $r_2(\overline {W(TR,1,J)}) =2n-1$.
\smallskip

\noindent (ii) By (5.7),  $r_2(\overline {W(TR,1,J)}) =2n-1$.
 This completes the proof of Lemma 5.2.\qed

\medskip
\noindent  Similar to Propositions 2.6 and  2.9, we 
 may prove that 

\smallskip
\noindent 
 {\bf Proposition 5.3.}
 {\em Suppose
 that $n\ge 4$ is even. 
 Then $r_2(\overline {S(T,R, J)}) = 2n-2$.
In particular, the rank of  $P\Omega
 _{2n}(\eta, q)$, where $n\ge 4$ is even,
  $\eta = \pm 1$, $q^n \equiv \eta\,\,
 (mod\,\, 4)$ is $2n-2$.
}
\smallskip

\noindent {\em Proof.} Since $S(T,R,J)$ is a Sylow
 2 subgroup of 
 $\Omega
 _{2n}(\eta, q)$ (where $n$ is even,
 $\eta = \pm 1$, $q^n \equiv \eta$
 (mod 4)), we conclude that 
$r_2(S(T,R ,J)) = 2n-1$ (see [L]).
 This implies that  $r_2(\overline {S(T,R ,J)}) 
\ge
r_2({S(T,R ,J)}) -1=
 2n-2.$ By Lemma 5.2,
 $$2n-1
=r_2(\overline {W(TR,1,J)})
  \ge r_2(\overline {S(T,R ,J)}) 
\ge 2n-2.$$ Since $T$ and $R$ (for $P\Omega_{2n}(\eta, q)$)
 satisfy the assumption of 
 Lemma 2.5, we have
 $r_2(\overline {S(T,R ,J)}) 
= 2n-2.$ \qed

\section{Orthogonal Commutator Groups $P\Omega
 _{2n}(\eta, q)$ where $n$ is odd,
 $\eta = \pm 1$, $q^n \equiv \eta$
 (mod 4)}
\noindent
 Let $n = 1+n_1$, where $n_1$ is even.
 Since  $P\Omega
 _{6}(\eta, q)$ ($n_1=2$) is not  a simple group, 
 we shall assume that $$n_1 \ge 4.$$
Let $S = S(T,R, J)$ be a Sylow 2-subgroup of 
$\Omega
 _{2n_1}(\eta', q)$ where $n_1$ is even,
 $\eta' = \pm 1$, $q^{n_1} \equiv \eta'$
 (mod 4) and let 
 diag$\,(e,1,1,\cdots , 1) = d(e)$,
 diag$\,(f,1,1,\cdots , 1)= d(f)$
 ($S$, $e$ and $f$ are given in section 5).
Let 
$$Y = \left < S, d(e), d(f)\right >
 \times \left < x, y\right >,
\mbox{ where } o(x)=o(y) = 2, o(xy) = 2^t,$$
and $ 2^{t+1}$ is the largest power of 2 that 
 divides $(q^2-1).$
Let $z = \mbox{diag}\,
(g^{2^{t-1}},g^{2^{t-1}},
 \cdots )$ be given as in section 5 (see (5.1)).
By Theorem 12 of [W],
 $ V =  \left < S, d(e)x,d(f)y\right >$
 is a Sylow 2-subgroup of 
$\Omega
 _{2n}(\eta, q)$
 and
$$
 \overline V = \left < S, d(e)x,d(f)y\right >
\Big /
 \left < z (xy)^{2^{t-1}}\right >\eqno(6.1)$$
 is a Sylow 2-subgroup of 
 $P\Omega
 _{2n}(\eta, q)$ ($n$ is odd,
 $\eta = \pm 1$, $q^n \equiv \eta$
 (mod 4)).
 It is easy to see that 
 $$\overline V \subseteq \overline W = 
\left <
W(TR,1,J), x, y \right >\Big /
\left < z (xy)^{2^{t-1}}\right >.\eqno(6.2)
$$
Note that $[W(TR,1,J) , \left <x,y\right> ]=1$. 
 As a consequence, 
$$\left <
W(TR,1,J), x, y \right >
\cong
\left (\begin{array}
{cc}
 W(TR,1,J)&0\\
0&1\\
\end{array}\right ) 
\times
\left  (\begin{array}
{cc}
 1&0\\
0&\left <g, k\right > \\
\end{array}\right ).\eqno(6.3)$$
Note that $\left <g, k\right >$ is 
 dihedral of order $2^{t+1}$ (see section 5) and that  
 $\left <x,y\right >\cong
 \left <g,h\right >$. It follows that 
one may decompose the group in (6.3) into
 $$(6.3)=
\left (\begin{array}
{cc}
 W(TR,1,J)&0\\
0&1\\
\end{array}\right ) 
\rtimes D = \Omega \rtimes D,\eqno(6.4)$$
 where $\Omega \cong  W(TR,1,J)$,
$D = \left
 <
\mbox{diag}\,
(v,v,\cdots)\,:\, v \in \left
<g, k\right > \right >.$
Note that $D \cong \left <g, k\right >$
 and that under the isomorphism we mention in
 (6.3), $z(xy)^{2^{t-1}}$ is mapped to 
$$ \sigma =
 \mbox{diag}\,( z,g^{2^{t-1}})=
 \mbox{diag}\,( g^{2^{t-1}}, g^{2^{t-1}},
 \cdots).\eqno(6.5)$$
Consequently,  
$$\overline W = 
  \overline {\Omega}\rtimes
\overline D
 \cong \Omega \rtimes
\overline D,
\eqno(6.6)$$
where  $\overline D = D/\left < \sigma\right >.$
 Since  $\Omega \cong W(TR,1,J)\cong \prod w_{m_i}
(TR)$ is a direct product of 
 wreath products ((4.8) of [L]),
we may apply Proposition 3.4 of [L] and conclude that 
$$r_2(\overline W)
 \le r_2(W(TR,1,J)) + r_2(D/\left < \sigma\right >)
 = 2n_1 +2.\eqno (6.7)$$
 Suppose that the above is actually an equality.
 Let $\overline E \subseteq \overline 
 \Omega $  be elementary abelian of dimension
 $2n_1 +2$.
Let $\overline E = (\overline E \cap \overline
 \Omega ) \times \overline F$.
 Applying 
 Lemma 2.3 of [L],  $\overline E \cap \overline
 {\Omega}$
 is of dimension $2n_1$
 and $r_2(\overline F)
 =2$. Applying Lemmas 2.3,  3.2 of [L],
  $\overline E \cap \overline {\Omega }$ is of the 
 form
$$\mbox{diag}\,
 (E_1, E_2, \cdots,
 E_{n_1/2} , 1 )\left <\sigma\right >/\left <\sigma\right >.
 \eqno(6.8)$$
 where $E_i \subseteq TR$ is of dimension 4 ($r_2(TR)=4$)
 and $r_2(\overline E|\overline
 D) = r_2(\overline F|\overline D) = r_2(\overline D)
 = 2$ (Lemma 2.3).  Hence  $\overline E|\overline D$
contains the centre of $\overline D$. 
Hence  $\overline E|\overline D$
must contain  the  element $\overline \tau$,
 $$ \tau =\mbox{diag}\,(v,v,\cdots, v),$$
$ v = g^{ \pm2^{t-2}}, \,\tau ^2 = \sigma$.
 It follows that the element in $\overline
 E$ projects to $\overline \tau$
 is of the form
$$\overline \nu
 =\overline {\mbox{diag}\,(  A  ,
    g^{ \pm2^{t-2}}          ) },\eqno (6.9) $$
 where $A \in W(TR, 1,J)$, $A^2 = z$
 ($A$ must satisfy the equation $A^2 = z$ as 
 $\overline \nu \in \overline E$ is of order 2).
Note that 
$[(6.8), (6.9)] = \overline 1.$
An easy observation
 of the  last entry of (6.8) and (6.9)
 shows that 
diag$\,(A, 1)\in C_{\Omega}(\mbox{diag}\,
(E_1, \cdots , E_{n_1/2},1 ))$
 (recall that $\overline \Omega \cong
 \Omega$).
 By the remark of Lemma 1.5, we have 
 $A \in 
\mbox{diag}\,(TR, TR, \cdots , TR)$.
 Since $A^2 = z$, we have
$$A=\mbox{diag}\,(a_1,a_2, \cdots),$$
 where $a_i $ is one of the members in (A2)
 (see Appendix A).
 Since $[a_i, E_i]=1$, the $a_i$ cannot be
$g^{2^{t-2}}d^mh,\,\,\,
d^{2^{t-2}}g^nk$ (see the discussion of 
 (5.5) and (5.6)
 of section 5). Hence 
 $a_i\in \{ g^{2^{t-2}},\,\,\, d^{2^{t-2}}\}.$
As a consequence,
 $$E_i \subseteq C_{TR}( g^{ \pm2^{t-2}} )
 \mbox{ or } C_{TR}( d^{ \pm2^{t-2}} ).$$
By Appendix A, such centralisers are 
 of rank 3.
This is a contradiction
 ($r_2(E_i) = 4 )$.
 Hence (6.7) can be refined into 
$$r_2(\overline W)
 \le 2n_1 +1.\eqno (6.10)$$
Applying (6.3), it is clear that 
$ r_2(\overline W) \ge r_2(W(TR,1,J)
 + 2-1 = 2n_1 +1$. Hence
$$r_2(\overline W)
  =  2n_1 +1.\eqno (6.11)$$
We are now ready to determine the rank of $\overline
V$.
 Note first that $V$ ($V$ is a 
 Sylow 2-subgroup of $\Omega_{2n_1}(\eta ', q)$)
  is of rank $2n_1 +1$.
 (see section 5.7 of [L]). It follows that 
$r_2(\overline V ) \ge 2n_1 +1 -1 = 2n_1.$
 By (6.11), we have 
$$2n_1 \le r(\overline V) \le 2n_1 + 1
 = r_2(\overline W).\eqno (6.12)$$
Suppose that  $r_2(\overline V) = 2n_1 + 1$.
 It follows that if $\overline E \subseteq
 \overline V$ is of maximal dimension, then 
 $\overline E$ is actually of maximal 
 dimension  in $\overline W$ ($r_2(\overline V) = 2n_1 + 1
 = r_2(\overline W))$.
 Since $r_2(\overline W) = 2n_1+1,
r_2(\overline {W(TR,1,J)}) = 2n_1-1$
 ((6.12) and Lemma 5.2), we have 
   $r_2(\overline E|\overline D) =2$ or 1. Suppose 
 that  $r_2(\overline E|\overline D) = 2$.
 \begin{enumerate}
\item[(i)] Similar to (6.9) , $\overline E$ possesses
 an element $\overline \tau,$ 
 $$\tau = \mbox{diag}\,(  A  ,
    g^{ \pm2^{t-2}}          ),\,\,
A \in W(TR, 1,J),\,\,A^2 = z.\eqno(6.13)$$
\item[(ii)]
 $\overline E \cap
 \overline \Omega $ is of rank $2n_1-1$.
\end{enumerate}
By the remark of Lemma 3.2 of 
 [L], 
 $\overline E \cap
 \overline \Omega  \subseteq
 \mbox{diag}\,(
 TR, TR, \cdots, TR, 1)\left <\sigma\right >/ 
\left <\sigma\right >  .$
 As a consequence, 
 $$\overline E \cap
 \overline \Omega =
 \mbox{diag}\,(
 E_1, E_2, \cdots, E_{n_1/2}, 1)
\left <\sigma\right >/ 
\left <\sigma\right >     ,\eqno(6.14)$$
 where all the $E_i$ are of rank
 4 except one which is of rank 3.
Note that 
$[(6.13), (6.14)] $
 $ =\overline 1.$
 It follows that 
$[A,  \mbox{diag}\,(
 E_1, E_2, \cdots, E_{n_1/2})]=1$.
By  the remark of Lemma 1.5, we have   
$ A \subseteq \mbox{diag}\,(TR, TR, \cdots , TR)$.
 Since $A^2 = \sigma$, we have
$$A=\mbox{diag}\,(a_1,a_2, \cdots),$$
 where $a_i $ is one of the members in (A2)
 (see Appendix A).
 Since $[a_i, E_i]=1$, the $a_i$ cannot be
$g^{2^{t-2}}d^mh,\,\,\,
d^{2^{t-2}}g^nk$  if $r_2(E_i) = 4$
(see the discussion of (5.5) and (5.6)
 of section 5). Hence 
 $a_i\in \{ g^{2^{t-2}},\,\,\, d^{2^{t-2}}\}$
 for all $i$ except one.
 Since $n_1 \ge 4$, at least one of the 
 $E_i$'s is of dimension  4. For such $E_i$,
 one has
 $$E_i \subseteq C_{TR}( g^{ \pm2^{t-2}} )
 \mbox{ or } C_{TR}( d^{ \pm2^{t-2}} ).$$
By Appendix A, such centralisers are 
 of rank 3.
This is a contradiction
 ($r_2(E_i) = 4 )$.
Hence $r_2(\overline E|\overline D) =1$
 and 
 $\overline E\cap \overline \Omega $ is 
 of dimension $2n_1$. In particular,
 $\overline E\cap \,\overline \Omega$ and  $\overline
\Omega $
 have the same dimension (recall that $\Omega \cong \overline 
 \Omega$). 
 By Lemma 3.2 of [L],
 $\overline E\cap \overline \Omega
 =  \mbox{diag}\,(
 E_1, E_2, \cdots, E_{n_1/2}, 1) \left <\sigma
 \right >/ \left <\sigma
 \right >  $, where
  $E_i \subseteq TR$ is  elementary abelian of 
 rank 4 for every $i$.
 In particular,
$\overline E$
 possesses an element of the 
 form
$$ \overline
 {\mbox{diag}\,(tr,1, \cdots, 1)} \in
\overline { S(T,R,J)}.$$
 Since $r_2(TR) =4$, $r_2(T)=3$, we may choose 
$r\in R^{\times}$.
 This is a contradiction
 (see remark of section 4.2 of [L]
 and the last paragraph of 
 the proof of Lemma 2.5). In summary, we 
 have
$$r_2(\overline V) = 2n_1.$$
 Equivalently, we have

\smallskip
\noindent 
 {\bf Proposition 6.1.}
 {\em The $2$-rank of 
 $P\Omega
 _{2n}(\eta, q)$ where $n\ge 5$ is odd,
 $\eta = \pm 1$, $q^n \equiv \eta$
 $($mod $4)$ is $2n_1 = 2(n-1)$.}

\bigskip
\begin{center}
{\bf Appendix A}
\end{center}

\medskip
The main purpose of this appendix is to 
 investigate the solutions of 
   the following  equation in $TR$
 (see section 5).
$$x^2 = z_0,\eqno(A1)$$
 where $Z(TR)
 = \left < z_0=  g^{2^{t-1}}\right >$.
 Note that $o(z_0) = 2$ and that 
 elements in $TR$ admit the following (unique)
 form
$g^n(dg)^mh^ok^p e^q f^r.$
 Direct calculation shows that 
$x$ is one of the following :
$$g^{2^{t-2}},\,\,\, d^{2^{t-2}} ,\,\,\,
g^{2^{t-2}}d^mh,\,\,\,
d^{2^{t-2}}g^nk,\,\, \mbox{ where } n, m \in \Bbb Z.
 \eqno(A2)$$
 
\smallskip
Let $w$ be an element in (A2) and let $C
 = C_{TR}(w)$ be its
 centraliser in $TR$. In the case $w= g^{2^{t-2}}$,
 one has 
 $$C = \left <
g, d, h,  kef\right >
= [(\left <g\right > \times 
 \left < dg \right >) \rtimes \left < h \right >]
\rtimes  \left < kef \right >.
$$
Suppose that the rank of $C$ is 4. Let
 $E \subseteq C$ be elementary abelian of 
 dimension  4.
Since $E$ is of dimension  4,
 $E \cap 
(\left <g\right > \times 
 \left < dg \right >)\rtimes
\left < h \right >
$ is  of dimension 3 (see Lemma 2.3 of [L]).
Hence  $E \cap 
(\left <g\right > \times 
 \left < dg \right >)$ must be of dimension 2
 (see Lemma 2.3 of [L]).
 Note that $\left <g\right > \times 
 \left < dg \right >$ has a unique
 elementary abelian subgroup of 
 dimension 2 : $\left <
 (dg)^{2^{t-2}}, (d^{-1}g)^{2^{t-2}}
 \right >$.
It follows that  $$\left <
 (dg)^{2^{t-2}}, (d^{-1}g)^{2^{t-2}}
 \right > \subseteq E$$ and that 
 $E$ must possess an element
 of the form $d^mg^n h$
 ( $r_2(E \cap 
(\left <g\right > \times 
 \left < dg \right >)\rtimes
\left < h \right >)=3>2=
r_2(E \cap 
(\left <g\right > \times 
 \left < dg \right >)$).
This is not possible as 
$d^mg^n h$ does not commute with
 $\left <
 (dg)^{2^{t-2}}, (d^{-1}g)^{2^{t-2}}
 \right >$.
 Hence the rank of $C$ is at most 3.
 One can show similarly that
 the rank of $C_{TR}(g^{2^{t-2}})$
 is
 at most 3 as well.

\bigskip
\bigskip
{\small
\noindent DEPARTMENT OF MATHEMATICS,\\
NATIONAL UNIVERSITY OF SINGAPORE,\\
SINGAPORE 117543,\\
REPUBLIC OF SINGAPORE.}

\bigskip

\noindent {\tt e-mail: matlml@math.nus.edu.sg}

\bigskip
\bigskip

\noindent  lang-68-4.tex

\end{document}